\documentclass[oneside]{amsart}
\usepackage[osf,sc]{mathpazo}
\usepackage[letterpaper, body={12.6cm, 19.0cm}, mag=1000]{geometry}
\usepackage{amsmath, amsthm, amssymb, enumerate, tikz, hyperref}
\usepackage{amsfonts}

\title{Finite groups with abelian automorphism groups :  A Survey}
\author{Rahul Dattatraya Kitture}
\author{Manoj K. Yadav}
\newtheorem{proposition}{Proposition}
\newtheorem{theorem}[proposition]{Theorem}
\newtheorem{remark}[proposition]{Remark}
\newtheorem{cor}[proposition]{Corollary}

\newcommand{\C}{\operatorname{C} }

\newcommand{\Z}{\operatorname{Z} }

\newcommand{\gen}[1]{\left < #1 \right >}
\newcommand{\Aut}{\operatorname{Aut} }

\newcommand{\Hom}{\operatorname{Hom} }

\newcommand{\Inn}{\operatorname{Inn} }

\newcommand{\Gcd}{\operatorname{gcd} }
\newcommand{\Autcent}{\operatorname{Autcent} }

\begin{document}

\subjclass[2010]{20D45, 20D15}
\keywords{Miller group, central automorphism, special $p$-group}

\maketitle
\begin{abstract}
In this  article we  present an extensive survey on the developments in the
theory of non-abelian finite groups with abelian automorphism groups, and pose
some problems and further research directions.
\end{abstract}
\section{Introduction}

All the groups considered here are finite groups.

Any cyclic group has abelian automorphism group. By the structure theorem for
finite abelian groups, it is easy to see that among abelian groups, only the
cyclic ones have abelian automorphism groups. A natural question, posed by H.
Hilton \cite[Appendix, Question 7]{hilton} in 1908, is:

\vspace{.1in}
{\it Can a non-abelian group have  abelian group of  automorphisms?}
\vspace{.1in}

An affirmative answer to this question was given by G. A. Miller \cite{miller}
in 1913.  He constructed a non-abelian group of order $2^6$ whose automorphism
group is elementary abelian of order $2^7$. Observe that a non-abelian group
with abelian automorphism group must be a  nilpotent group of nilpotency class
$2$. Hence it suffices to study the groups of prime power orders.
Investigation on the structure of groups with abelian automorphism groups was
initiated by C. Hopkins \cite{hopkins} in 1927. He, among other things, proved
that such a group can not have a non-trivial abelian direct factor, and if such
a group is a $p$-group, then so is its  automorphism group, where $p$ is a prime
integer.

Unfortunately, the topic was not investigated for about half a century after the
work of Hopkins. But days came when examples of such odd prime power order
groups were constructed by D. Jonah and M. Konvisser \cite{jonah} in 1975 and a
thesis was written on the topic by B. E. Earnley \cite{earnley} during the same
year, in which, attributing to G. A. Miller, a group with abelian automorphism
group was named as `Miller group'. Following Earnley, we call a group to be {\it
Miller} if it is non-abelian and its automorphism group is abelian. Earnley
pointed out that the former statement of Hopkins is not true as such for
$2$-groups. He constructed Miller $2$-groups admitting a non-trivial abelian
direct factor. However, the statement is correct for odd order groups. He also
presented a generalization of examples of Jonah-Konvisser, which implies that
the number of elements in a minimal generating set for a Miller group can be arbitrarily
large. But he obtained a lower bound on the number of generators, which
was later improved to an optimal bound by M. Morigi \cite{morigi2}. He also gave
a lower bound on the order of Miller $p$-groups, which was again improved to an
optimal bound by M. Morigi \cite{morigi1}.

Motivated by the work done on the topic, various examples of Miller groups were
constructed by several mathematicians via different approaches during the next
20 years, which were mostly special $p$-groups. A $p$-group $G$ is said to be
{\it special} if $\Z(G) = G'$ is elementary abelian, where $\Z(G)$ and $G'$
denote the center and the commutator subgroup of $G$ respectively. These include
the works by R. Faudree \cite{faudree1}, H. Heineken and M. Liebeck
\cite{heineken2}, D. Jonah and M. Konvisser \cite{jonah}, B. E. Earnley
\cite{earnley}, H. Heineken \cite{heineken1}, A. Hughes \cite{hughes}, R. R.
Struik \cite{struik}, S. P. Glasby \cite{glasby}, M. J. Curran \cite{curran1},
and M. Morigi \cite{morigi1}. The existence of non-special Miller groups follows
from \cite[Remark 2]{jonah}.

The abelian $p$-groups of minimum order, which can occur as automorphism groups
of some $p$-groups were studied by P. V. Hegarty \cite{hegarty} in 1995 and G.
Ban and S. Yu \cite{ban} in 1998. More examples of Miller groups were
constructed by A. Jamali \cite{jamali} and Curran \cite{curran3}.

Neglecting Remark 2 in \cite{jonah}, A. Mahalanobis \cite{mahalanobis}, while
studying Miller groups in the context of MOR cryptosystems, conjectured that
Miller $p$-groups are all special for odd $p$. Again  neglecting \cite[Remark
2]{jonah}, non-special Miller $p$-groups were constructed by V. K. Jain and the
second author \cite{yadav1}, V. K. Jain, P. K. Rai and the second author
\cite{yadav2}, A. Caranti \cite{caranti} and the authors \cite{yadav3}.

The motivation for many of the examples comes from some natural questions or
observations on previously known examples. The construction of examples of
Miller groups of varied nature have greatly contributed to complexify  the
structure of a Miller group. The structure of a general Miller group  has not
yet been well understood.

We record here the known information about
structure of Miller $p$-groups and their automorphism groups.
Let $G$ be a Miller $p$-group, and assume that it has no abelian direct factor.
Then the following hold:
\begin{enumerate}
\item $|G|\geq p^6$ if $p=2$ and $|G|\geq p^7$ if $p>2$.
\item Minimal number of generators of $G$ is $3$ for $p$ even, and $4$ for
$p$ odd.
\item The exponent of $G$ is at least  $p^2$.
\item If $|G'|>2$, then $G'$ has at least two cyclic factors of the maximum
order in its cyclic  decomposition.
\item If ${\rm Aut}(G)$ is elementary abelian, then $\Phi(G)$ is elementary
abelian, $G'\leq \Phi(G)\leq Z(G)$ and at least one equality holds. There are
Miller $p$-groups in which only one equality holds.
\item For $p>2$ (resp. $p=2$), the abelian $p$-groups of order $<p^{11}$ (resp.
$<2^7$) are not automorphism groups of any $p$-group.
\end{enumerate}

The aim of this article is to present an extensive survey on the developments in
the theory of Miller groups since 1908, and pose some problems and further
research directions. There are, at least, three survey articles on automorphisms
of $p$-groups and related topics (see \cite{helleloid}, \cite{malinowska},
\cite{mann}); the present one does not overlap with any of them.


We conclude  this section with  setting some notations for a multiplicatively
written group $G$.  We denote by $\Phi(G)$,  the Frattini subgroup of  $G$. For
an element $x$ of $G$, $\langle x\rangle$ denotes the cyclic subgroup generated
by $x$, and $o(x)$ denotes its order. For subgroups $H,K$ of a group $G$, $H<K$
(or $K>H$) denotes that $H$ is proper subgroup of $K$. The exponent of   $G$ is
denoted by $\exp(G)$. By $\C_n$, we denote the cyclic group of order $n$. For a
$p$-group $G$ and integer $i\geq 1$, $\Omega_i(G)$ denotes the subgroup of $G$
generated by those $x$ in $G$ with $x^{p^i}=1$, and $\mho^i(G)$ denotes the
subgroup generated by $x^{p^i}$ for all $x$ in $G$. By ${\Inn}(G)$,  ${\Autcent}(G)$ and ${\Aut}(G)$ 
we denote, respectively, the group of inner automorphisms,
central automorphism, and all automorphisms of $G$. Let a group $K$ acts on a
group $H$ by automorphisms, then $H \rtimes K$ denotes the semidirect product of
$H$ by $K$. All other notations are standard.
%
%
\section{Reductions}

Starting from the fundamental observations of Hopkins, in this section, we
summarize all the known results (to the best of our knowledge) describing the
structure of Miller groups.
Although we preserve the meaning,  no efforts are made to preserve the original
statements from the source.

As commented in the introduction,  a non-abelian Miller group must be nilpotent
(of class $2$), and therefore it is sufficient to  study Miller  $p$-groups. The
following result is an easy exercise.

\begin{proposition}[Hopkins, \cite{hopkins}]
Every automorphism of a Miller group centralizes $G'$.
\end{proposition}

An automorphism of a group $G$ is said to be {\it central} if it induces the
identity automorphism on the central quotient $G/\Z(G)$. Note that
${\Autcent}(G)$, the
group of all  central automorphisms of $G$, is the centralizer of ${\Inn}(G)$
 in ${\Aut}(G)$.

\vskip3mm
 A group $G$ is said to be {\it purely
non-abelian} if it has no non-trivial abelian direct factor. Hopkins made the
following
important observation:
\begin{theorem}[Hopkins, \cite{hopkins}]\label{hopkins1}
 {\it A Miller $p$-group is purely non-abelian.}
\end{theorem}
Unfortunately, the statement, as such, is not true for $p=2$ as shown by Earnley
\cite{earnley} (see Theorem \ref{earnley1} below for the correct statement).

\vskip3mm
In a finite abelian group $G$, there exists a minimal generating set
$\{x_1,\ldots,x_n\}$ such that $\langle x_i\rangle \cap \langle
x_j\rangle =1$ for $i\neq j$. Hopkins observed that Miller groups possess a set
of generators with a property similar to  the preceding one. More precisely, he
proved the following result.

\begin{theorem}[Hopkins, \cite{hopkins}]
If $G$ is a Miller $p$-group, then there exists a set $\{x_1,\ldots,x_n\}$ of
generators of $G$ such that
\begin{enumerate}
\item for $p>2$, $\langle x_i\rangle \cap \langle x_j\rangle=1$ for all
$i\neq j$;
\item for $p=2$, $\langle x_i\rangle \cap \langle x_j\rangle$ is of order at
most $2$, for all $i\neq j$.
\end{enumerate}
\end{theorem}
In fact the theorem is true for {\it any $p$-group of class $2$}, as shown
below. Let $G$ be a $p$-group of class $2$. Choose a set $\{y_1,\ldots, y_k\}$
of generators of $G$ with the property

\vspace{.1in}
{\it $\prod_{i} o(y_i)$ is minimum.} \hfill (*)
\vspace{.1in}\\
Let $p>2$  and  $o(y_i)=p^{n_i}$, $1 \le i \le k$. If $\langle y_1\rangle \cap
\langle
y_2\rangle\neq 1$ then $y_2^{p^{n_2-1}}\in \langle y_1\rangle$. Assuming
that $o(y_1)\geq o(y_2)$, we can write $y_2^{p^{n_2-1}}=(y_1^{p^{n_2-1}})^a$
for some integer $a$. Take $y_2'=y_2y_1^{-a}$. Then $\{y_1,y_2',y_3,\ldots,
y_k\}$ is a generating set for $G$ and $o(y_2')\leq p^{n_2-1}<o(y_2)$, which
 contradicts  (*). The case $p=2$ can be handled in a similar way.

\vspace{.1in}
\begin{theorem}[Hopkins, \cite{hopkins}]
If $G$ is a Miller $p$-group, then  ${\Aut}(G)$ is a $p$-group.
\end{theorem}

This can be proved easily in the following way. In the case when $G$ is
purely non-abelian, by  a result of Adney-Yen \cite[Theorem 1]{adney},
${\Autcent}(G)$
($={\Aut}(G)$) has order equal to $|{\Hom}(G,\Z(G))|$,
which is clearly a power of $p$. Now consider the case when $G$ has an
abelian direct factor, which occurs only when  $p=2$ and it must be cyclic of
order at least $2^2$ (see Theorem \ref{earnley1}).
 Let $G=H\times C_{2^n}$, where $H$ is purely non-abelian and $n \ge 2$. By the
main theorem
in \cite{bidwell},
$$|{\Aut}(G)|=|{\Aut}(H)|\cdot |{\Aut}(C_{2^n})|\cdot |{\Hom}(H,C_{2^n})| \cdot
|{\Hom}(C_{2^n}, Z(H))|,$$
and each factor on the right side has order a power of $2$.

\vspace{.1in}
The investigation of structure of Miller groups remained unattained for about
half a century until it was revisited by
 Earnley \cite{earnley} in 1974. He  pointed out  that a Miller
$2$-group {\it can} have an abelian direct factor. So the correct form of
Theorem \ref{hopkins1} is

\begin{theorem}[Earnley, \cite{earnley}] \label{earnley1}
Let $G$ be a finite $p$-group such that $G=A\times N$ with $A\neq 1$ an abelian
group and $N$ a purely non-abelian group. Then $G$ is a Miller group  if
and only if $p=2$ and $A,N$ satisfy the following conditions:
\begin{enumerate}
\item $A$ is cyclic of order $2^n>2$;
\item $N$ is a special Miller $2$-group.
\end{enumerate}
\end{theorem}

\begin{theorem}[Earnley, \cite{earnley}] \label{earnley2}
If $G$ is a Miller $p$-group, then the following hold.
\begin{enumerate}
\item The exponent of $G$  is greater than $p$;
\item If $p>2$, then $Z(G)\cap \Phi(G)$ is non-cyclic.
\end{enumerate}
\end{theorem}

For the first statement of the preceding theorem, we can assume that $p>2$ and
the proof now follows by noting that if $\exp(G) = p$, then the
map $x\mapsto x^{-1}$ is a non-central  automorphism of $G$. The second
statement follows from the
following result  of Adney-Yen.

\begin{theorem}[Adney-Yen, \cite{adney}]\label{adney1}
If $G$ is a $p$-group of class $2$ such that $G'$ has only one
cyclic factor of maximum order in the direct product decomposition, then $G$
possesses a non-central automorphism.
\end{theorem}

Since, in $p$-group of class $2$, $G'  \le \Z(G )\cap \Phi(G)$,  the preceding
theorem  restricts   the structure of   the commutator subgroup of a Miller
$p$-group.

\begin{theorem}
If $G$ is a Miller $p$-group, $p$ odd, then $G'$ possesses at least
two cyclic factors of maximum order in the direct product decomposition.
\end{theorem}

This raises a natural question for $p=2$. The analogue of the preceding
theorem for $p=2$ holds except when $|G'|=2$. This can be obtained
from the following generalization of Theorem \ref{adney1}.

\begin{theorem}[Faudree, \cite{faudree}]
Let $G$ be a $p$-group of class $2$ with the following conditions:
\begin{enumerate}
\item $G'=\langle u\rangle \times {\rm U}$, where
$o(u)=p^{m_1}>p^{m'}=exp({\rm U}).$
\item $[g,h]=u$ and $h^{p^{m_1+m'}}=1$,
\item $m''=m'$ if $p$ is odd, and $m''=\mbox{max}(1,m')$ if $p=2$.
\end{enumerate}
Let $H=\langle g,h\rangle$ and $L=\{ x\in G : [g,x], [h,x]\in {\rm U}\}$. Then
$G=HL$
and the map
$$g\mapsto gh^{p^{m''}}, \hskip5mm h\mapsto h, \hskip5mm x\mapsto x\;\; (x\in
L)$$
defines an automorphism of $G$ which centralizes  $Z(G)$.
\end{theorem}

Notice that if $p=2$ and the exponent of $G'$ is at least  $4$, then the
automorphism defined in the preceding theorem is  non-central, and therefore $G$
is not Miller. But if $|G'|= 2$, then the theorem is no-longer applicable to
produce a non-central automorphism of $G$.  So the following question remains
open.

\vspace{.1in}
\noindent{\bf Question 1.} {Can  a finite $2$-group $G$ with $G'$ cyclic of
order $2$ be
Miller?}

\vspace{.1in}
Earnley obtained lower bounds for the order and the minimum number of
generator of a Miller group, which were later sharpened to the optimal level by
Morigi \cite{morigi1, morigi2}

\begin{theorem}
Let $G$ be a Miller $p$-group.
\begin{enumerate}
\item For any prime $p$, $G$ is generated by at least $3$ elements.
(Earnley, \cite{earnley})
\item If $p$ is odd, then $G$ is generated by at least $4$ elements.
(Morigi,\cite{morigi1})
\end{enumerate}
\end{theorem}

The example of a Miller $2$-group constructed by Miller is minimally generated
by $3$ elements (see Section 3  (3.1)).  For $p$ odd, Jonah-Konvisser
constructed Miller $p$-groups
which are minimally generated by $4$ elements (see Section 4 (4.3) for more
details).

\vspace{.1in}
\begin{theorem}\label{theorem11}
Let $G$ be a Miller $p$-group.
\begin{enumerate}
\item For any prime $p$, $|G|\geq p^6$. (Earnley, \cite{earnley})
\item If $p$ is odd, then $|G|\geq p^7$ and there exists a Miller
group of order $p^7$. (Morigi, \cite{morigi2})
\end{enumerate}
\end{theorem}

Again the Miller $2$-group constructed by Miller is of order $2^6$
having automorphism group of order $2^7$. Morigi constructed a special Miller
$p$-group $G$ of order $p^7$ with $\Aut(G)$ elementary abelian of order
$p^{12}$. In fact, it is one among an infinite family of Miller groups
constructed by Morigi (see Section 4 (4.6) for more details).

The question whether an abelian $p$-group  of order smaller than $p^{12}$ for an
odd prime  $p$ can occur as the automorphism group of a $p$-group, was addressed
by  Hegarty \cite{hegarty} and  Ban and   Yu \cite{ban} independently.
They proved

\begin{theorem} For $p$ odd, there is no abelian $p$-group of order smaller than
$p^{12}$
which can occur as the  automorphism group of a $p$-group.
\end{theorem}

Finally we state some results on the structure of ${\Aut}(G)$ for a Miller
group $G$. Many known examples of Miller groups are special $p$-groups.
Certainly for such groups, the automorphism group is elementary abelian.
The following problem appears as an {\it old problem} in
\cite[problem 722]{berkovich}.

\vspace{.1in}
\noindent {\bf Problem 2.} Study the $p$-groups with elementary abelian
automorphism groups.
\vspace{.1in}

Let $G$ be a $p$-group. If ${\Aut}(G)$ is elementary abelian, then so is
$G/Z(G)$; hence $G'\leq \Phi(G)\leq Z(G)$. In fact, it is interesting to see
that at least one equality always holds.

\begin{theorem}[Jain-Rai-Yadav, \cite{yadav2}]\label{jain-rai-yadav}
Let $G$ be a $p$-group, $p$ odd, such that ${\Aut}(G)$ is elementary abelian.
Then $\Phi(G)$ is elementary abelian, and one of the following holds:
\begin{enumerate}
\item $Z(G)=\Phi(G)$.
\item $G'=\Phi(G)$.
\end{enumerate}
\end{theorem}
Jain-Rai-Yadav \cite{yadav2} constructed $p$-groups with elementary abelian
automorphism group, in which exactly one of the above two conditions holds.

\vspace{.1in}
For $2$-groups with elementary abelian automorphism groups, there are analogous
necessary conditions, stated below. Since a Miller $2$-group can have abelian
(cyclic) direct factor, we consider two cases.

An abelian $p$-group of type $(p^n,p,\ldots,p)$ with $n>1$ is called a {\it
ce-group}. If $G=A\times B$ with $A \cong \C_{p^n}$ ($n>1$) and $B\cong
\C_p\times
\cdots \times \C_p$, then call $A$ {\it cyclic part} and $B$ {\it elementary
part} of
$G$.

\begin{theorem}[Jafari, \cite{jafari}]
Let $G$ be a purely non-abelian $2$-group. Then ${\Autcent}(G)$ is
elementary abelian if and only if one of the following holds.
\begin{enumerate}
\item $G/G'$ is of exponent $2$.
\item $Z(G)$ is of exponent $2$.
\item ${\Gcd}\big(\exp(G/G'), \exp( Z(G))\big)=4$ and $G/G', Z(G)$ are
ce-groups such that elementary part of $Z(G)$ is contained in $G'$ and there
is an element $z$ of order $4$  in cyclic part of $Z(G)$ with $zG'$ lying in
cyclic part of $G/G'$ satisfying $o(zG')=\exp(G/G')/2$.
\end{enumerate}
\end{theorem}
In particular, if $G$ is a purely non-abelian $2$-group with ${\Aut}(G)$
elementary abelian, then one of the above three conditions holds. The example
constructed
by Miller shows that it satisfies only condition (3). Jain-Rai-Yadav
\cite{yadav2}
constructed purely non-abelian Miller $2$-groups which satisfy only
condition (1) or only condition (2). Finally, for $2$-groups $G$ with abelian
direct factor, the following theorem gives necessary conditions for ${\rm
Aut}(G)$ to be elementary abelian.

\begin{theorem}[Karimi-Farimani, \cite{karimi}]
 Let $G=A\times N$ with $A$ cyclic
$2$-group and $N$ purely non-abelian $2$-group of class $2$. Then ${\rm
Autcent}(G)$ is elementary abelian if and only if
\begin{enumerate}
\item $|A|$ is $4$ or $8$.
\item $N$ is special $2$-group.
\end{enumerate}
\end{theorem}
With the notations of the preceding theorem, if $G$ is Miller then (1)
and (2) hold.

Theorem \ref{earnley1} tells us that a Miller $p$-group, $p$ odd,  can not admit
a non-abelian direct factor. It is natural to ask whether a Miller $p$-group can
occur as a direct product of non-abelian groups. This situation has been
considered by Curran \cite{curran3}. He determined necessary and
sufficient conditions on a direct product $H\times K$ to have abelian
automorphism group. Note that for any non-trivial group $H$, $\Aut(H \times H)$
is non-abelian. Before we state the result of Curran, we set some notations.

Let $H$ be a $p$-group of class $2$.

\begin{enumerate}
\item Let $a,b,c,d$ denote the exponents of $H/H'$, $H/\Z(H)$, $\Z(H)$ and
$H'$ respectively.
\item If $H'$ and $\Z(H)$ have the same rank,
then define $d_s$ to be the largest integer \big($\le \exp(\Z(H)$\big) such that
$\Omega_{d_s}(H')=\Omega_{d_s}(Z(H))$.
\item If $H/\Z(H)$ and $H/H'$ have the same rank, then define $b_t$ to be
the largest integer \big($\le \exp(H/H')$\big) such that
$ \Omega_{b_t}(H/\Z(H)) \cong\Omega_{b_t}(H/H')$.
\item Let $\Gamma^i(H)$ denote the subgroup of $H$ with
$\Gamma^i(H)/H'=\mho^i(H/H')$.
\end{enumerate}
Replacing $H$ by $K$ in (1)-(4), the corresponding terms $a',b',c',d'$,
$d'_{s'}$, $b'_{t'}$ and $\Gamma^i(K)$ are similarly defined. For simplicity,
we denote by $r(A)$, the rank of an abelian group $A$.  With this setting, we
have

\begin{theorem}[Curran, \cite{curran3}]
Let $G=H \times K$, where $H,K$ are $p$-groups of
class $2$ with no common direct factor. Then ${\Aut}(G)$ is abelian if and
only if ${\Aut}(H)$ and ${\Aut}(K)$ are abelian and one of the following
holds:
\begin{enumerate}
\item $\Z(H)=H'$ and $\Z(K)=K'$.
\item $\Z(H)>H'$ and $\Z(K)=K'$, where $r(\Z(H))=r(H')$,
$a'\leq d_s\leq a\leq c$ and $\Omega_a(Z(H))\leq \Gamma^{c'}(H)$.
\item $\Z(H)=H'$ and $\Z(K)>K'$, where $r(\Z(K))=r(K')$, $a\leq d'_{s'}\leq
a'\leq c'$ and $\Omega_{a'}(\Z(K))\leq \Gamma^c(K)$.
\item $\Z(H)>H'$ and $\Z(K)>K'$, where $r(\Z(H))=r(H')$, $r(\Z(K))=r(K')$ and
$a=d_s=d'_{s'}=a'$.
\item $\Z(H)>H'$ and $\Z(K)=K'$, where $r(H/\Z(H))=r(H/H')$, $a'\leq b_t\leq
c\leq a$ and $\Omega_{a'}(\Z(H))\leq \Gamma^c(H)$.
\item $\Z(H)=H'$ and $\Z(K)>K'$, where $r(K/\Z(K))=r(K/K')$, $a\leq b'_{t'}\leq
c'\leq a'$ and $\Omega_a(\Z(K))\leq \Gamma^{c'}(K)$.
\item $\Z(H)>H'$ and $\Z(K)>K'$, where $r(H/\Z(H))=r(H/H')$,
$r(K/\Z(K))=r(K/K')$ and
$c=b_t=b'_{t'}=c'$.
\end{enumerate}
\end{theorem}
%
\begin{remark}{\rm Observe that in the above theorem, if ${\Aut}(H \times K)$
is abelian then either $H'$ and $\Z(H)$ have the same rank or
$H/H'$ and $H/Z(H)$ have the same rank; the same is true for the other component
$K$.}
\end{remark}

An analogous problem for central product may be stated as follows.

\vspace{.1in}
\noindent {\bf Problem 3.} Find necessary and/or sufficient condition such that
central product of two Miller groups is again a Miller group.

%

\section{Examples of Miller $2$-groups}

In this section, we discuss  examples of Miller $2$-groups in chronological
order.
\vspace{.1in}

\noindent{\bf (3.1)}\label{section3.1}
As mentioned above several times, the first example of a
Miller $2$-group was constructed by  Miller himself, which comes as a
semi-direct product of the
cyclic group of order $8$ by the dihedral group of order $8$, and presented by
\begin{align*}
 G_1&=C_8\rtimes D_8 = \langle x,y,z \mid x^8, y^4, z^2, zyz^{-1}=y^{-1},
yxy^{-1}=x^5, zxz^{-1}=x\rangle.
\end{align*}
Miller proved that each coset of $Z(G_1)$ is invariant under every automorphism
of $G_1$ and  that every automorphism has order dividing $2$. Thus ${\Aut}(G_1)$
is elementary abelian, and it can be shown that its order is $2^7$.

By Theorem \ref{theorem11}, the order of a Miller $2$-group is at least
$64$. Having an example of order $64$, a natural idea which peeps in ones mind
is to explore groups of order $64$ to find more Miller groups. This was done by
Earnley \cite{earnley}, who proved that there are (exactly) two more Miller
groups of the minimal order, which are presented as follows.
\begin{align*}
 G_2 &= (C_4 \rtimes C_4) \rtimes C_4 = \langle x,y,z \mid yxy^{-1}=x^{-1},
zxz^{-1}=xy^2, zyz^{-1}=y\rangle,\\
G_3 &= (C_4\times C_4\times C_2)\rtimes C_2  = \langle x,y,z,t \mid x^4, y^4,
z^2,
t^2, xy=yx, xz=zx, yz=zy, \\
 & \hskip6cm txt^{-1}=xy^2, tyt^{-1}=tz, t^{-2}zt=z\rangle.
\end{align*}

Note that in all the groups $G_1,G_2$ and $G_3$, the center and Frattini
subgroups coincide and are elementary abelian of order $8$, whereas the
commutator subgroup is elementary abelian of order $4$. Further,
${\Aut}(G_2)$ and ${\Aut}(G_3)$ are elementary abelian  of order $2^9$.

A generalization of $G_1$ appears, as an exercise,  in the book
\cite[Exercise 46, Page 237]{macdonald} by Macdonald, written in 1970. For
$n\geq 3$, the group is presented as follows.
\begin{align*}
G_{1,n} = \langle a,b,c \mid & a^{2^n}=b^4=c^2=1,
b^{-1}ab=a^{1+2^{n-1}},\hskip2mm
c^{-1}bc=b^{-1},\hskip2mm [c,a]=1\rangle.
\end{align*}
The group $G_{1,n}$ is of order $2^{n+3}$ and its automorphism group is
an abelian $2$-group of type $(2^{n-2},2,2,2,2,2,2)$. In 1982,  Struik
\cite{struik}, independently, obtained the same example with a different
presentation.

A generalization of $G_2$ and $G_3$ has also been obtained by Glasby
\cite{glasby}, which is described as follows. For $n\geq 3$, let $G_{2,n}$
denotes the $2$-group of class $2$ with generators
$x_1,x_2,\cdots, x_n$ with following additional relations:
\begin{align*}
x_i^4=1, \,\, (1\leq i\leq n), \,\,  [x_i,x_n]=x_{i+1}^2,\,\, (1\leq i\leq
n-1),
\end{align*}
and set $[x_k,x_l]=1$ in the remaining cases.
Here ${\rm Aut}(G_{2,n})$ is elementary abelian group of order $2^{2n}$. If
$n=3$ then  $G_{2,n}$ is isomorphic to $G_2$.

Again for $n\geq 3$, let $G_{3,n}$ denotes the $2$-group of class $2$ with
generators $y_0$, $y_1$, $\ldots$, $y_n$ with following additional relations:
\begin{align*}
y_0^2=y_i^4=y_n^2=1, \,\, (1\leq i\leq n-1), \,\, [y_i,y_n]=y_{i+1}^2,\,\,
(1\leq i\leq n-2)
\end{align*}
and set $[y_k,y_l]=1$ in the remaining cases.
Here also ${\rm Aut}(G_{3,n})$ is elementary abelian group of order
$2^{2n}$. If $n=3$ then  $G_{3,n}$ is isomorphic to $G_3$.

\vspace{.1in}
\noindent{\bf (3.2)}\label{section3.3}
It should be noted that, in 1974,  Jonah and Konvisser constructed special
Miller groups of order $p^8$, which were generalized to an infinite family of
Miller $p$-groups by Earnley, and those groups include the case $p=2$ too (see
Section 4 (4.3) and (4.4)).

\vspace{.1in}
\noindent{\bf (3.3)}\label{section3.4}
Heineken and Liebeck (see Section 4 (4.2) for details) proved that {\it given a finite group $K$, there
exists a special $p$-group $G$, $p$ odd, with ${\rm Aut}(G)/{\rm
Autcent}(G) \cong K.$} In particular, if $K=1$, then the corresponding group
$G$ is Miller $p$-group.

In 1980, A. Hughes \cite{hughes} proved that
one can construct a special $2$-group as well with the above property. The
method of Hughes is a modification of the graph theoretic method of
Heineken-Liebeck. It should be noted that the method of
Heineken-Liebeck  uses {\it digraphs}, whereas that of  Hughes considers {\it graphs}, and
is described below.

Let $K$ be any finite group and associate to $K$ a connected graph $D(K)$ which
satisfies following conditions:
\begin{enumerate}
\item Each vertex of the graph has degree at least $2$.
\item Every cycle in the graph contains at least $4$ vertices.
\item ${\rm Aut}(D(K))\cong K$.
\end{enumerate}
Associate a special $2$-group $G$ to $D(K)$ as follows. If the graph has $n$
vertices $v_1,\ldots,v_n$, then consider the free group $F_n$ on
$x_1,\ldots,x_n$. Let $R$ be the normal subgroup of $F_n$ generated by $x_i^2$,
$[x_i,[x_j,x_k]]$ (for all $i,j,k$) and $[x_r,x_s]$ whenever the vertices $v_r$ and
$v_s$ are adjacent. Define $G$ to be the group $F_n/R$; it is a special
$2$-group of order $2^{n+\binom{n}{2}-e}$, where $n$ is the number of generators
of $G$ (so $|G/G'|=2^n$) and $e$ is the number of edges of the graph.  It turns out that ${\rm Aut}(G)/{\rm
Autcent}(G)\cong K$ (see \cite{hughes} for proof).

\vspace{.1in}
\noindent{\bf (3.4)}\label{section3.2} In 1987,  Curran \cite{curran1} studied
automorphisms of
semi-direct product,
and suggested a method to construct many more examples of Miller $2$-groups
similar to $G_1$. We describe the method briefly and see that the above $G(n)$
can be constructed by this method. Let $A=\langle a\rangle$ be a cyclic group of
order $2^n$, $n\geq 3$ and  $N$  a special $2$-group acting on $A$ in  the
following way: a maximal subgroup $J\leq N$ acts trivially, and any $x\in
N\setminus J$ acts by
$$xax^{-1}=a^{1+2^{n-1}}.$$
Let $G = A\rtimes N$ be the semi-direct product of $A$ and $N$ with this action.
Then we get

\begin{theorem}[Curran, \cite{curran1}]\label{curthm}
Let $G = A\rtimes N$ be as above along with the conditions:
\begin{enumerate}
\item $A\times J$ is characteristic in $A\rtimes N$.
\item Any automorphism of $N$ leaving $J$ invariant is central automorphism
of $N$.
\end{enumerate}
Then $A\rtimes N$ is a Miller group.
\end{theorem}

To elaborate, consider $A=\langle a\rangle$, the cyclic group of order
$2^n$, $n\geq 3$ and $N=\langle b,c \mid b^4,c^2,[b,c]=b^2\rangle$, the dihedral
group of order $8$. Consider the action of $N$ on $A$ by
$$b^{-1}ab=a^{1+2^{n-1}}, \hskip5mm c^{-1}ac=a.$$
Define $G=A\rtimes N$, the semi-direct product with this
action. Note that $Z(G)=\langle a^2,b^2\rangle=\Phi(G)$, and $G/\Phi(G)$ is
elementary abelian of order $8$.

Take $J=\langle b^2,c\rangle$, the largest subgroup of $N$ acting trivially
on $A$. Then $A\times J$ is characteristic in $A\rtimes N$, since it is the
unique abelian subgroup of index $2$ (if there were more than one abelian
subgroups of index $2$, then center would have index $4$). Also, in $N$, the
subgroup $\langle b\rangle$ is characteristic.
Hence if $\varphi\in {\Aut}(N)$ leaves $J$ invariant, then it leaves
invariant the subgroup $J\cap
\langle b\rangle=\langle b^2\rangle$ (the center of $N$), and one can see
that $\varphi$ is central automorphism of $N$. The group $A\rtimes N$ is
therefore a Miller group by Theorem \ref{curthm}.
Note that this group is isomorphic to $G(n)$ described above.

\vspace{.1in}
\noindent{\bf (3.5)}\label{section3.5}
After a considerable time gap, Jamali \cite{jamali} in 2002 constructed the
following infinite family of Miller $2$-groups in which,
the number of generators and the exponent of the group can be arbitrarily large.
For integers $m\geq 2$ and $n\geq 3$, let $G_n(m)$ be the group generated by
$a_1,\ldots, a_n,b$, subject to following relations:
\begin{align*}
&a_1^2=a_2^{2^m}=a_i^4=1 \hskip5mm (3\leq i\leq n)\\
 & a_{n-1}^2=b^2,\\
& [a_1,b]=[a_i,a_j]=1 \;\;\;(1\leq i<j\leq n),\\
& [a_n,b]=a_1, [a_{i-1},b]=a_i^2 \;\;\;(3\leq i\leq n).
\end{align*}
The group $G_n(m)$ has order $2^{2n+m-2}$ and its automorphism
group is  abelian of type $(\underbrace{2,2,\ldots,2}_{n^2}, 2^{m-2})$.
For $n=3$ and $m=2$, the group $G_3(2)$ is isomorphic to $G_3$ (see Section
3(3.1)).

\section{Examples of Miller $p$-groups,  $p$-odd}

This section, a lifeline for Miller groups in a sense, presents  evolution of
the topic. We'll see the  influence of examples of varied nature  on Miller
groups towards understanding the structure.
Some examples of Miller groups occurred in other related contexts  without any
pointer to the topic.
The deriving force behind many of the  examples comes from  natural
questions on previously known Miller groups or certain natural optimistic
expectations on the structure
of such groups.
The evolution process certainly helped, although minimally,  in  studying the
structure of Miller
groups, especially, in turning down the natural optimistic expectations.

During 1971 to 1979, there were three occasions, in which certain $p$-groups
were constructed with a specific property, and these groups turned out to
be Miller groups.  The motive of the
 construction of these groups had no obvious connection with Miller groups.

\vspace{.1in}
\noindent{\bf (4.1)}\label{section4.1}
We now start  the discussion of the (possibly) first example of Miller group  of
odd order. It was conjectured that a finite group in which every element
commutes with its epimorphic image is abelian. It was disproved 
by R. Faudree \cite{faudree1} in 1971. He constructed  a {\it non-abelian} group $G$ in
which every element commutes with its epimorphic image.  The group $G$ is a
special $p$-group, and therefore one can easily
deduce that ${\rm Aut}(G)$ is abelian. The group $G$ is described as follows.

Let $G=\langle a_1,a_2,a_3,a_4\rangle$ be the $p$-group of class $2$ with
following additional relations:
\begin{align*}
& [a_1,a_2]=a_1^p, \, [a_1,a_3]=a_3^p, \, [a_1,a_4]=a_4^p,\\
& [a_2,a_3]=a_2^p, \, [a_2,a_4]=1, \, [a_3,a_4]=a_3^p.
\end{align*}
The group $G$ is a special $p$-group of order $p^8$ and $Z(G)=G'$ is elementary
abelian of order $p^4$. {\it If $p$ is odd}, then ${\rm Aut}(G)$ is elementary
abelian of order $p^{16}$.

\vspace{.1in}
\noindent{\bf (4.2)}\label{section4.2}
In 1974,  Heineken-Liebeck \cite{heineken2},  constructed a $p$-group $G$ of class $2$, $p$ odd,
for a given finite group $K$ such that ${\Aut}(G)/{\Autcent}(G) \cong K$.
Associate to the finite group $K$, a connected {\it digraph} (directed graph)
$X$ as follows. If $K$ is cyclic  of order $>2$, then take $X$ to be the cyclic
digraph with $|K|$ vertices. In the remaining cases, we associate a digraph $X$
with following conditions:

(1) $X$ is {\it strongly} connected,

(2) Any two non-simple vertices are {\it not} adjacent,

(3) Every vertex belongs to a (directed) cycle of length at least $5$,

(4) Every non-simple vertex has at least two outgoing edge,

(5) ${\rm Aut}(X)\cong K$,\\
where by a {\it simple} vertex we mean a vertex with exactly one incoming and
exactly one outgoing edge. Associate to such $X$, a $p$-group $G_X$ of class $2$
as follows. If $X$ has $n$ vertices, we take $n$ generators for $G_X$. If the
vertex $i$ has outgoing {\it edges}  to $j_1,\ldots,j_r$ precisely,  then we put
the relation $x_i^p=[x_i, x_{j_1}\cdots x_{j_r}]$. Finally we put
$[[x_i,x_j],x_k]=1$ for all $i,j,k$, making $G_X$ of class $2$. It is easy to
show that $G_X$ is a special $p$-group of order $p^{n+\binom{n}{2}}$.

\vspace{.1in}
The conditions (1)-(4) imply that ${\rm Aut}(G)/{\rm Autcent}(G)$ is
isomorphic to the automorphism group of $X$, which is isomorphic to $K$ (by
(5)).

\vspace{.1in}
Before proceeding for the example of Miller $p$-group by this method, we make
some comments. If $|K|\geq 5$, then there is
a systematic procedure to construct graph $X$ satisfying conditions (1) to
(5); it is obtained by a specific subdivision of a Cayley digraph of $K$.
If $|K|<5$, then  one constructs $X$ satisfying (1)-(5) by some ad-hoc method.

\vspace{.1in}
For $K=1$, the following digraph satisfies conditions (1)-(5).
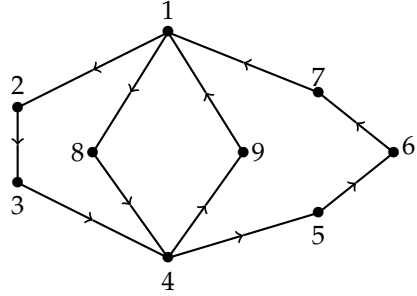
\begin{center}
\begin{tikzpicture}
\draw [->, thick] (0,0) -- (1,0.3);\draw [thick] (1,0.3) -- (2,0.6);
\draw [thick, ->] (2,0.6) -- (2.5,1);\draw [thick] (2.5,1) -- (3,1.4);
\draw [thick, ->] (3,1.4) -- (2.5,1.8);\draw [thick] (2.5,1.8) -- (2,2.2);
\draw [thick, ->] (2,2.2) -- (1,2.6);\draw [thick] (1,2.6) -- (0,3);
\draw [->, thick] (0,3) -- (-1,2.5);\draw [thick] (-1,2.5) -- (-2,2);
\draw [thick, ->] (-2,2) -- (-2,1.5);\draw [thick] (-2,1.5) -- (-2,1);
\draw [thick, ->] (-2,1) -- (-1,0.5);\draw [thick] (-1,0.5) -- (0,0);
\draw [->, thick](0,0) -- (0.5,0.7);\draw [thick] (0.5,0.7) -- (1,1.4);
\draw [->, thick](1,1.4) -- (0.5,2.2);\draw [thick] (0.5,2.2) -- (0,3);
\draw [thick,->] (0,3) -- (-0.5,2.2);\draw [thick] (-0.5,2.2) -- (-1,1.4);
\draw [thick, ->] (-1,1.4) -- (-0.5,0.7); \draw [thick] (-0.5,0.7) -- (0,0);
\node at (1,1.4) {$\bullet$}; \node at (-1,1.4) {$\bullet$};
\node at (0,0) {$\bullet$};\node at (2,0.6) {$\bullet$};
\node at (3,1.4) {$\bullet$};
\node at (2,2.2) {$\bullet$};
\node at (0,3) {$\bullet$};
\node at (-2,1) {$\bullet$};
\node at (-2,2) {$\bullet$};
\node at (0,3.3) {$1$};
\node at (-2,2.3) {$2$};
\node at (-2,0.7) {$3$};
\node at (0,-0.3) {$4$};
\node at (2,0.3) {$5$};
\node at (3.2,1.5) {$6$};
\node at (2,2.4) {$7$};
\node at (-1.2, 1.4) {$8$};
\node at (1.2, 1.4) {$9$};
\end{tikzpicture}
\end{center}

The group $G_X$ associated to this digraph is a special Miller $p$-group of
order $p^{9+\binom{9}{2}}=p^{45}$, $p$ odd.

\vspace{.1in}
\noindent{\bf (4.3)}\label{section4.3}
The existence of Miller $2$-groups of order $2^6$ served as a base to the
belief that $p^{45}$ is too big to be a minimal order for a Miller $p$-group for
odd $p$. This motivated  Jonah and  Konvisser \cite{jonah} to construct
Miller $p$-groups of smaller orders.  In 1975, they constructed, for {\it each
prime} $p$, $p+1$ non-isomorphic Miller $p$-groups of order $p^8$ as described
below.   Let $\lambda=(\lambda_1,\lambda_2)$ be a {\it non-zero} vector with
entries in field $\mathbb{F}_p$ and $G_{\lambda}=\langle a_1,a_2,b_1,b_2\rangle$
be the $p$-group of class $2$ with following additional relations:
\begin{align*}
&a_1^p=[a_1,b_1], \; a_2^p=[a_1,b_1^{\lambda_1}b_2^{\lambda_2}],\;
b_1^p=[a_2,b_1b_2],\\
& b_2^p=[a_2,b_2], \; [a_1,a_2]=[b_1,b_2]=1.
\end{align*}
The proof of the fact that $G_{\lambda}$ is Miller group is very elegant.  We
briefly
describe the idea here as  it will  be useful in later discussions.

Fix a non-zero vector $\lambda$ and write $G=G_{\lambda}$. Note that $G$
is special $p$-group of order $p^8$ and $G'$ is of rank $4$ with generators
$[a_i,b_j]$, $i,j=1,2$.

The subgroups $A=\langle a_1,a_2,G'\rangle$ and $B=\langle b_1,b_2,G'\rangle$
are the only abelian normal subgroups such that $[A:G']= p^2 = [B:G']$;
hence they are permuted by every automorphism of $G$. There exists $x\in A$ such
that $A^p\leq [x, G]$, but there is no $y\in B$ with $B^p\leq [y,G]$; hence
$A,B$ are characteristic in $G$. The only $x\in A$ with $A^p\leq [x,G]$ are
$a_1^iz$ with $p\nmid i$ and $z\in G'$; they generate characteristic subgroup
$\langle a_1,G'\rangle$. Similarly, those $x\in A$ with $B^p\leq [x,G]$ generate
a characteristic subgroup, namely $\langle a_2,G'\rangle$ and those
$x\in B$ with $x^p\in [x,G]$ generate the characteristic subgroup $\langle
b_2,G'\rangle$. Thus, if $\varphi\in{\rm Aut}(G)$, then
\begin{equation*}
\varphi(a_1)\equiv a_1^i, \hskip5mm
\varphi(a_2)\equiv a_2^j, \hskip5mm
\varphi(b_1)\equiv b_1^{k_1}b_2^{k_2}, \hskip5mm
\varphi(b_2)\equiv b_2^l   \pmod {G'},
\end{equation*}
where $i,j,l$ are integers not divisible by $p$, and $(k_1,k_2)\neq (0,0)$.
Using relations in $G$, it is now easy to deduce that $i=j=k_1=l=1$ mod $p$ and
$k_2=0$. Therefore every automorphism of
$G$ is central. Since $G'=Z(G)$, it follows that ${\rm Aut}(G)$ is abelian.

The groups $G_{\lambda}$ and
$G_{\mu}$, for non-zero vectors $\lambda,\mu\in\mathbb{F}_p\times
\mathbb{F}_p$, are isomorphic only if these vectors are dependent.

\begin{remark}\label{rem1}
Note that in the group $G_{\lambda}$, replacing, in the power relations, all the
$p$-th power by $p^k$-th power, it can be shown by the same arguments that the
resulting groups are still Miller. If $k>1$, then
$G_{\lambda}'=Z(G_{\lambda})$ and has exponent $p^k>p$; hence $G_{\lambda}$ is a
{\it{non-special Miller group}}.
\end{remark}

We must emphasize the negligence of the preceding remark  \cite[Remark 2]{jonah}
by the authors of many recent papers claiming the existence of non-special
Miller groups was not known in the literature until recently.

\vspace{.1in}
\noindent{\bf (4.4)}\label{section4.4}
Examples of Jonah and Konvisser were generalized by  Earnley
\cite{earnley}, during the same year, by increasing
number of generators, in the following way.
Fix $n\geq 2$ and a non-zero $n$-tuple $(\lambda_1,\ldots,\lambda_n)$ with
$\lambda_i\in\mathbb{F}_p$. Let $G_{\lambda}=\langle x_1,x_2, y_1,y_2,\ldots,
y_n\rangle$ be the $p$-group of class $2$ with the following additional
relations:
\begin{align*}
& x_1^p=[x_1,y_1], \hskip5mm x_2^p=[x_1,y_1^{\lambda_1}\ldots
y_n^{\lambda_n}],\\
& y_i^p=[x_2, y_iy_{i+1}], \hskip3mm i=1,2,\ldots, n-1,\\
& y_n^p=[x_2,y_n],
 [x_1,x_2]=[y_i,y_j]=1 \hskip3mm \mbox{ for all }i,j.
\end{align*}
It can be shown easily that $G_{\lambda}'=Z(G_{\lambda})$ is elementary
abelian and  generated by the $2n$ elements
$[x_i,y_j]$ for $i=1,2$ and $j=1,2,\ldots,n$; hence order of $G_{\lambda}$ is
$p^{(n+2)+2n}=p^{2+3n}$.

For ${\Aut}(G_{\lambda})$ to be abelian, the obvious necessary condition is
that {\it every automorphism should be central}, but since $G_{\lambda}$ is
special, this is sufficient too. We briefly describe the beautiful linear
algebra techniques evolved by Earnley to prove that every automorphism of
$G_{\lambda}$ is central as, with a little variation,  these techniques have
been used in different contexts by Morigi \cite{morigi1}, Hegarty \cite{hegarty}
and Earnley himself.

For simplicity, fix non-zero vector $\lambda$ and write $G=G_{\lambda}$. Let
$f:G/G'\rightarrow G'$ denote the map $xG'\mapsto x^p$, which
 is a homomorphism for odd $p$. An automorphism $\alpha$ of $G/G'$ determines
its action on $G'$ by
$$\hat{\alpha}:G'\rightarrow G', \hskip5mm \hat{\alpha}([x,y])=[\alpha(xG'),
\alpha(yG')].$$
The automorphism $\alpha$ of $G/G'$ is induced by an automorphism $\varphi$ of
$G$ if and only if $\hat{\alpha}\circ f = f\circ\alpha$, i.e. the following
diagram commutes:
\begin{center}
\begin{tikzpicture}
\draw [->, thick] (1,0) -- (2,0);
\draw [->, thick] (1,2) -- (2,2);
\draw [<-, thick] (0,0.5) -- (0,1.5);
\draw [<-, thick] (3,0.5) -- (3,1.5);
\node at (0,0) {$G/G'$};\node at (0,2) {$G/G'$};
\node at (3,0) {$G'$.};\node at (3,2) {$G'$};
\node at (1.5,0.3) {$f$};\node at (1.5,2.3) {$f$};
\node at (0.5,1) {$\alpha$};\node at (3.5,1) {$\hat{\alpha}$};
\end{tikzpicture}
\end{center}
Consequently, every automorphism of $G$ is central if and only if  identity is
the only automorphism of $G/G'$ which fits in the above commutative diagram.

Note that $G/G'$ and $G'$ can be considered as vector spaces over
$\mathbb{F}_p$; hence the automorphisms $\hat{\alpha}$ and $\alpha$ in the above
diagram can be considered as (invertible) linear maps. However, if $p=2$, the
map $f$ {\it may not} be linear. But, in the group $G$ under consideration,
consider the {\it abelian} subgroups $A=\langle x_1,x_2,G'\rangle$ and
$B=\langle y_1,\ldots, y_n,G'\rangle$; they generate $G$ and the restriction of
$f$ to $A/G'$ and $B/G'$ are homomorphisms, and therefore linear.

Next, $A$ and $B$ are the only abelian subgroups containing $G'$ such that
modulo $G'$ each of them have order at least $p^2$. For $n>2$, both $A$ and $B$
are characteristic
as their orders are different. But even for $n=2$, they are characteristic,
since, in this case, $A$ contains an element $t$ such that $A^p\leq [t,G]$, but
$B$ has no element with similar property.

Thus, consider an automorphism $\alpha=(\alpha_1,\alpha_2)$ of $G/G'=A/G'
\oplus B/G'$, where $\alpha_1 \in {\Aut}(A/G')$ and $\alpha_2\in{\Aut}(B/G')$.
As a vector space, $G'=[A,B]$ is isomorphic to $A/G' \otimes B/G'$, the tensor
product of $A/G' $ and $B/G'$. Hence the automorphism induced by
$\alpha$ on $G'$ is nothing but $\alpha_1\otimes \alpha_2$. Then $\alpha$ is
 induced by an automorphism of $G$ if and only if the following diagrams
commute:

\begin{center}
\begin{tikzpicture}
\node at (-0.5,0) {$A/G'$};\node at (-0.5,2) {$A/G'$};
\node at (2.5,0) {$G'$};\node at (2.5,2) {$G'$};
\draw [->, thick] (0,0) -- (2,0);
\draw [->, thick] (0,2) -- (2,2);
\draw [->, thick] (-0.5,1.5) -- (-0.5,0.5);
\draw [->, thick] (2.5,1.5) -- (2.5,0.5);
\node at (1,0.5) {$f$};\node at (1,2.5) {$f$};
\node at (-0.7,1) {$\alpha_1$};\node at (3.2,1) {$\alpha_1\otimes\alpha_2$};
\node at (5.5,0) {$B/G'$};\node at (5.5,2) {$B/G'$};
\node at (8.5,0) {$G'$};\node at (8.5,2) {$G'$};
\draw [->, thick] (6,0) -- (8,0);
\draw [->, thick] (6,2) -- (8,2);
\draw [->, thick] (5.5,1.5) -- (5.5,0.5);
\draw [->, thick] (8.5,1.5) -- (8.5,0.5);
\node at (7,0.5) {$f$};\node at (7,2.5) {$f$};
\node at (5.3,1) {$\alpha_2$};\node at (9.2,1) {$\alpha_1\otimes\alpha_2$};
\end{tikzpicture}
\end{center}
For the first diagram, consider the following ordered bases:
\begin{align*}
\{ x_1G', x_2G'\} \mbox{ for } A/G'  \hskip3mm\mbox{ and }\hskip3mm
\{[x_1,y_1], [x_2,y_1], \ldots, [x_1,y_n], [x_2,y_n]\}
\mbox{ for } G'.
\end{align*}
Also for the second diagram, consider the following bases:
\begin{align*}
\{ y_1G', \ldots, y_nG'\} \mbox{ for } B/G'   \hskip2mm\mbox{ and }\hskip2mm
\{[x_1,y_1], \ldots, [x_1,y_n], [x_2,y_1],\ldots,
[x_2,y_n]\} \mbox{ for } G'.
\end{align*}

The matrix of $f$ with respect to these bases can be easily written from the
{\it power-commutator relations} in $G$. Writing the matrices of
$\alpha_1,\alpha_2,\alpha_1\otimes\alpha_2$ with respect to these bases, one
expresses the commutativity of above diagram in terms of two matrix equations,
and a simple matrix computation shows that $\alpha_1 = 1$ and
$\alpha_2 =1$ are the only solutions.
This implies that every automorphism of $G$ is central.
Note that this also covers the case $p=2$ as  the restrictions of $f$ to $A/G'$
and $B/G'$  are  linear.

 For more detailed module-theoretic formulation of the above arguments for
arbitrary special $p$-groups, one may refer to \cite{caranti}.

\vspace{.1in}
\noindent{\bf (4.5)}\label{section4.5}
The third instance of examples of Miller groups, not in the context of the
topic, is  by  Heineken \cite{heineken1} in 1979. He constructed a family of
$p$-groups in which every normal subgroup is characteristic. The groups,
actually, possess more interesting properties, which force  the groups actually
become Miller. The  construction is as follows. Let $\mathbb{F}_q$ ($q=p^n$) be
the field of order $p^n$ and ${\rm U}(3,\mathbb{F}_q)$ denote the group of
$3\times  3$ unitriangular matrices over $\mathbb{F}_q$. Identify the elements
of ${\rm U}(3,\mathbb{F}_q)$ as triples over $\mathbb{F}_q$ by
$$(x,y,z) \leftrightsquigarrow
\begin{bmatrix}
1 & x & z\\
 & 1 & y\\
  & & 1
\end{bmatrix}, \hskip5mm x,y,z\in\mathbb{F}_q.
$$

The sets $A=\{ (x,0,z): x,z\in
\mathbb{F}_q\}$ and $B=\{ (0,y,z): y,z\in \mathbb{F}_q\}$ constitute abelian
subgroups of order $q^2=p^{2n}$ and exponent $p$; they generate ${\rm
U}(3,\mathbb{F}_q)$ and their intersection is the center of
${\rm U}(3,\mathbb{F}_q)$. Also the commutator of an element of $A$ with
that of $B$ is given by $$[(x,0,z_1), (0,y,z_2)]=(0,0,xy).$$

Heineken constructed a group $G$ by a little modification in the commutator
and power relations of ${\rm U}(3,\mathbb{F}_q)$ as follows. The group $G$
consists of triples $(x,y,z)$ over $\mathbb{F}_q$, in
which $$A^*=\{ (x,0,z) \mid x,z\in\mathbb{F}_{q}\} \mbox{ and }
B^*=\{ (0,y,z) \mid y,z\in\mathbb{F}_{q}\} $$
constitute abelian subgroups of order $q^{2}=p^{2n}$. The power relations are
\begin{align*}
(x,0,z)^p &= (0,0,x^p-x),\\
(0,y,z)^p &= (0,0,y+y^p+\cdots + y^{p^{n-1}}).
\end{align*}
For a fixed generator $t$ of the multiplicative group of $\mathbb{F}_{q}$, the
commutator relations between $A^*$ and $B^*$ are defined as
$$[(x,0,z_1), (0,y,z_2)]=(0,0,xy-tx^py^{p^2}).$$

Then $G=\{(x,y,z) \mid x,y,z\in\mathbb{F}_{q}\}$ is a group, which is the
product
of abelian subgroups $A^*$ and $B^*$, and $Z(G)=\{ (0,0,z) \mid
z\in\mathbb{F}_{q}\} =
A^*\cap B^*.$
Fix $x\neq 0$ and vary $y\in\mathbb{F}_{q}$, then the elements
$xy-tx^py^{p^2}$ exhaust whole $\mathbb{F}_{q}$. It follows that
$G'=Z(G)$, and therefore  $G$ is a special $p$-group of
order $p^{3n}$. With the above setup, we have

\begin{theorem}
Let $n$ be odd positive integer. If $n\geq 5$, $p>2$  or $n=3$, $p\geq
5$, then every automorphism of $G$ is identity on $Z(G)$ as
well as on $G/Z(G)$; hence ${\Aut}(G)$ is abelian.
\end{theorem}

The groups $G$ in the preceding theorem have the following remarkable property:
for every $x \in G-G'$, the conjugacy class of $x$ in $G$ is $xG'$. The groups
satisfying this property are called {\it Camina} groups. The readers interested
in Camina groups are referred to \cite{DS96} and the references therein.

An automorphism $\alpha$ of a group is said to be {\it class-preserving} if it
takes each element of the group to its conjugate.
It is an easy exercise to show that each automorphism of the group $G$ is
class-preserving.
This is  not only true for the groups $G$
considered above, but also for any Camina $p$-group which is Miller as well.

\begin{theorem}
Let $G$ be a  $p$-group which is  Miller as well as a Camina group. Then
automorphisms of $G$ are all class-preserving.
\end{theorem}

As a simple consequence, we have
\begin{cor}
Let $G$ be a  $p$-group which is  Miller as well as a Camina group. Then normal
subgroups of $G$ are all characteristic.
\end{cor}

We remark that the examples of Jonah-Konvisser are not Camina groups and
extraspecial $p$-groups are Camina groups but not Miller. As the  structure of
Camina  $p$-groups of class $2$ in general is not well understood, it will be
interesting  to study
$p$-groups which are both Camina  as well as Miller.

\vspace{.1in}
\noindent{\bf Problem 4.} Determine the structure of  $p$-groups which are
Camina
 as well as Miller.

\vspace{.1in}
This structural information will also, on the one hand, help in understanding
$p$-groups
of class $2$ whose automorphisms  are all class-preserving and, on the other
hand, shed some light on the study of $p$-groups of class $2$ in which all
normal subgroups are characteristic.

\vspace{.1in}
\noindent{\bf (4.6)}\label{section4.6}
Note that the Miller groups constructed by Heineken in the preceding discussion
are of order at least $p^9$. But, as we already mentioned, the lower bound on
the order of Miller groups is $p^7$ for odd $p$. This was  Morigi
\cite{morigi1}, who, in 1994, constructed examples of Miller groups of minimal
order as a part of a general construction of an infinite family of such groups.
The construction is briefly described as follows.
For any natural number $n$, let $G^n$ denote the $p$-group of class $2$
generated by $a_1,a_2,b_1,\ldots, b_{2n}$ with the following additional
commutator and power relations:
\begin{align*}
& [a_1, b_{2i+1}]=[a_2, b_{2i+2}]=1, \hskip5mm i=0,1,\ldots, n-1,\\
& [b_1,b_2] = [b_3,b_4] = \cdots = [b_{2n-1}, b_{2n}]=1,\\
& [b_i,b_j]=1 \mbox{ if } i\equiv j \pmod 2.
\end{align*}
Further $a_1^p=a_2^p=1$ and $b_1^p$ is the product of following $n^2+n+1$
commutators:
\begin{equation*}
\tag{**} [a_1, a_2], [a_1, b_{2i+2}], [a_2, b_{2i+1}], [b_{2i+1}, b_{2j+2}],
\hskip5mm i,j=0,1,\ldots, n-1, i\neq j.
\end{equation*}
Finally powers of $b_i$'s are related by
\begin{align*}
& b_2^p = b_1^p [ a_1,b_2]^{-1}\\
& b_{2i+1}^p = b_{2i}^p [a_2, b_{2i+1}]^{-1}, \,\,\,\, b_{2i+2}^p =
b_{2i+1}^p[a_1, b_{2i+2}]^{-1}, \hskip5mm i=1,\ldots, n-1.
\end{align*}
Then $G^n$ is a special $p$-group, with $|G^n/(G^n)'|=p^{2n+2}$ and $(G^n)'$ is
elementary
abelian $p$-group generated by  $n^2+n+1$ commutators in $(**)$; so
$|G^n|=p^{n^2+3n+3}$. Further, ${\Aut}(G^n)$ is (elementary) abelian, and
$$ |{\Aut}(G^n) |= |{\Autcent}(G^n)|= p^{(2n+2)(n^2+n+1)}.$$

For $n=1$, $G^1$ is a special $p$-group of order $p^7$.  This is an example of a
Miller group of the smallest order for
$p$ odd. The following problem is interesting:

\vspace{.1in}
\noindent {\bf Problem 5.} Describe all Miller $p$-groups of order $p^7$ for  $p
>
2$.

\vspace{.1in}
\noindent{\bf (4.7)}\label{section4.7}
Upto this point, we have only considered the examples of special Miller
$p$-groups (modulo Remark \ref{rem1}). Now we'll consider non-special groups.
Before we proceed further with more examples, we record a very useful result of
 Adney and  Yen \cite{adney}

Let $G$ be a purely non-abelian $p$-group of class $2$. Let $G/G'=\prod_{i=1}^r
\langle x_iG'\rangle$, with $o(x_iG')\geq
o(x_{i+1}G')$. Let $p^a, p^b, p^c$ denote the exponents of $Z(G)$, $G'$ and
$G/G'$ respectively. Finally, let $R$ be the subgroup of $Z(G)$ generated by
all the homomorphic images of $G$ in $Z(G)$, and $K$ denote the intersection of
the kernels of all the homomorphisms $G\rightarrow G'$.

\begin{theorem}[Adney-Yen, \cite{adney}]\label{adney-yen}
With the above notations, ${\Autcent}(G)$ is abelian if and only if the
following hold:
\begin{enumerate}
\item $R=K$.
\item either min$(a,c)=b$ or min$(a,c)>b$ and $R/G'=\langle
x^{p^b}G'\rangle$.
\end{enumerate}
\end{theorem}

Motivated by the conjecture of Mahalanobis as stated in the introduction,
Jain and the second author \cite{yadav1}, in 2012,  constructed the
following infinite family of non-special Miller
$p$-groups. For $n\geq 2$, and $p$ odd, let $G_n=\langle
x_1,x_2,x_3,x_4\rangle$ be the $p$-group of class $2$ with the following
additional relations:
\begin{align*}
& x_1^{p^n}=x_2^{p^2}=x_3^{p^2}=x_4^p=1,\\
& [x_1,x_2]=x_2^p, \; [x_1,x_3]=[x_1,x_4]=x_3^p,\\
& [x_2,x_3]=x_1^{p^{n-1}}, \; [x_2,x_4]=x_2^p,\;
 [x_3,x_4]=1.
\end{align*}
It is then easy to see that
\begin{enumerate}
\item $Z(G_n)=\Phi(G_n)=\langle x_1^p, x_2^p, x_3^p\rangle$.
\item $G_n'=\langle x_1^{p^{n-1}}, x_2^p, x_3^p\rangle$ is elementary
abelian of order $p^3$.
\item $G_n$ is special only when $n=2$ (follows by (1) and (2)).
\end{enumerate}
The proof that $\Aut(G_n) = \Autcent(G_n)$  is constructive and reply on
detailed careful
calculations. An application of Theorem \ref{adney-yen} now shows that $G_n$ is
Miller.

\vspace{.1in}

In the preceding examples of non-special Miller $p$-groups $G$, one can easily
notice that $$G' < \Z(G) = \Phi(G).$$
 One might desire that for a Miller $p$-group $G$, one of the following always
holds true: (i) $G' = \Z(G)$; (ii) $\Z(G) = \Phi(G)$.

In 2013, Jain, Rai and the second author \cite{yadav2} constructed the following
 infinite family of Miller $p$-groups $G$
such that $G'<Z(G)<\Phi(G)$, which we again denote by $G_n$. For $n \ge 4$, let
$G_n=\langle x_1,x_2,x_3,x_4\rangle$ be a $p$-group of class $2$ with the
following additional
relations:
\begin{align*}
& x_1^{p^n}=x_2^{p^4}=x_3^{p^4}=x_4^{p^2}=1\\
& [x_1,x_2]=[x_1,x_3]=x_2^{p^2}, \; [x_1,x_4]=x_3^{p^2},\\
& [x_2,x_3]=x_1^{p^{n-2}}, \; [x_2,x_4]=x_3^{p^2},\;
 [x_3,x_4]=x_2^{p^2}.
\end{align*}
Then $G_n$ is a $p$-group of order $p^{n+10}$ with
\begin{align*}
Z(G_n) =\langle x_1^{p^2}, x_2^{p^2}, x_3^{p^2}\rangle ,\;
\Phi(G_n) = \langle x_1^p, x_2^p, x_3^p, x_4^p\rangle ,\;
G_n' = \langle x_1^{p^{n-2}}, x_2^{p^2}, x_3^{p^2}\rangle.
\end{align*}
It follows that
$G_n'\leq Z(G_n)<\Phi(G_n)$ and $G_n'=Z(G_n)$ only when $n=4$.
As in the previous examples, the proof of the fact that $G_n$ is Miller is
constructive.

\vspace{.1in}
As mentioned in Theorem \ref{jain-rai-yadav}, if ${\rm
Aut}(G)$ is elementary abelian, then $\Phi(G)$ is elementary abelian and one of
the following holds: (1) $Z(G)=\Phi(G)$; (2) $G'=\Phi(G)$.
Jain, Rai and the second author constructed the following Miller $p$-groups in
which exactly one of these conditions hold.

\vspace{.1in}
For any prime $p$, let $G_4=\langle x_1,x_2,x_3,x_4\rangle$ denote the
$p$-group of class $2$ with the following additional relations:
\begin{align*}
& x_1^{p^2}=x_2^{p^2}=x_3^{p^2}=x_4^{p^2}=1,\\
& [x_1,x_2]=1, \, [x_1,x_3]=x_4^p, \, [x_1,x_4]=x_4^p,\\
& [x_2,x_3]=x_1^p,\,[x_2,x_4]=x_2^p,\, [x_3,x_4]=x_4^p.
\end{align*}
Then $G_4$ is a $p$-group of order $p^8$ and the following holds.
\begin{enumerate}
\item $G_4'=\langle x_1^p,x_2^p,x_4^p\rangle$ is elementary
abelian of order $p^3$.
\item $\Phi(G_4)=Z(G_4)=\langle x_1^p,x_2^p,x_3^p,x_4^p\rangle$ is
elementary abelian of order $p^4$.
\item ${\Aut}(G_4)$ is elementary abelian of order $p^{16}$.
\end{enumerate}

Again for any prime $p$, consider the $p$-group $G_5=\langle
x_1,x_2,x_3,x_4,x_5\rangle$ of class $2$ with the following additional
relations:
\begin{align*}
& x_1^{p^2}=x_2^{p^2}=x_3^{p^2}=x_4^{p^2}=x_5^p=1,\\
& [x_1,x_2]=x_1^p, \, [x_1,x_3]=x_3^p, \, [x_1,x_4]=1,\, [x_1,x_5]=x_1^p,\,
[x_2,x_3]=x_2^p,\\
&[x_2,x_4]=1, \, [x_2,x_5]=x_4^p, \, [x_3,x_4]=1, \, [x_3,x_5]=x_4^p,\,
[x_4,x_5]=1.
\end{align*}
Then $G_5$ is a $p$-group of order $p^9$ and the following holds.
\begin{enumerate}
\item $G_5'=\Phi(G_5)=\langle x_1^p,x_2^p,x_3^p,x_4^p\rangle$ is elementary
abelian of order $p^4$.
\item $Z(G_5)=\langle x_4,G_5'\rangle$.
\item ${\Aut}(G_5)$ is elementary abelian of order $p^{20}$.
\end{enumerate}

It is natural to ask

\noindent {\bf Question 6.}
Does there exist a Miller $p$-group $G$ in which $Z(G)\nsubseteq \Phi(G)$ and
$\Phi(G) \nsubseteq Z(G)$?

\vspace{.1in}
\noindent{\bf (4.8)}\label{section4.8}
The proofs of the results in subsection (4.7)  involve heavy computations. To
remedy the problem, in 2015, Caranti
\cite{caranti, caranti2} suggested
a simple module theoretic approach to construct non-special  Miller $p$-groups
from special ones. The arguments given in \cite{caranti} are not sufficient to
prove the results as stated. The authors of the present survey proved that the
results are valid under an additional hypothesis. The construction  is briefly
described as follows.

For an odd prime $p$, let $H$ be a special Miller $p$-group satisfying the
following hypotheses:
\begin{enumerate}
\item[(i)] $\mho_1(H)$ is a proper subgroup of $H'$.
\item[(ii)] The map $H/H'\rightarrow H'$ defined by $hH'\mapsto h^p$ is
injective.
\end{enumerate}
Let $K = \gen{z}$ be the cyclic group of order $p^2$, and $M$ be a subgroup of
order $p$ in $H'$ but not in $\mho_1(H)$.  Let $G$ be a central product of $H$
and $K$ amalgamated at $M$. Note that $G'=\Phi(G)<Z(G)$; hence $G$ is
non-special. With this setting,  we have

\begin{theorem}
 If $H/M$ is a Miller group, then so is $G$.
\end{theorem}

\vspace{.1in}
Before we proceed, we make a comment on the preceding theorem. Caranti  claimed
that $G$ is Miller without the condition `$H/M$ is Miller'. Unfortunately, this
is not always true, as shown in the  following example.

\vspace{.1in}
Let $H=\langle a,b,c,d\rangle$ be the $p$-group of class
$2$ with the following additional relations:
\begin{equation*}
a^p=[a,c], \hskip3mm b^p=[a,bcd], \hskip3mm c^p=[b,cd],\hskip3mm d^p=[b,d].
\end{equation*}
Then  $H$ is a special Miller $p$-group of order $p^{10}$ and satisfies
conditions (i)-(ii).
It can be proved  that if $M = \gen{[a,b]}$,
then $G$ is Miller, and  if $M = \gen{[a,d]}$, then $G$ is not a Miller group.
That $H$ and $H/\gen{[a,b]}$ are Miller can be proved
 following the arguments similar to those in  \cite{jonah} (for details see
\cite{yadav3}).

\vspace{.1in}
As noted above, the  Miller groups $G$  are such that $G'=\Phi(G)<Z(G)$. Now
with  a little variation in the preceding construction, we obtain Miller
$p$-groups $G$ with $G'<\Phi(G)=Z(G)$. Assume that the special Miller group $H$
also
satisfies the following hypothesis in addition to (i)-(ii) above:
\begin{enumerate}
\item[(iii)] If $H$ is minimally generated by $x_1,\ldots,x_n$, then $H'$ is
minimally generated by $[x_i,x_j]$ for $1\leq i<j\leq n$.
\end{enumerate}
Let $L = \gen{z}$ be a cyclic group of order $p^n$, $n \ge 3$ and let $z$ act on
$H$ via a non-inner central
automorphism $\sigma$ of $H$ (which always exists in a Miller $p$-group). Let
$N$
be a subgroup of order $p$ in $H'$ but not in $\mho_1(H)$. Now define $G$ to be
the {\it partial semi-direct product} of $H$ by $L$ amalgamated at $N$ (cf.
\cite{gor} or \cite{yadav3}). With this setting, we finally have

\begin{theorem}
 If $H/N$ is a Miller group, then so is $G$.
\end{theorem}
\vspace{.1in}
Again, we remark that the condition `$H/N$ is Miller', in the preceding theorem,
can not be dropped,  as shown by the following example.
Let $H=\langle
a,b,c,d\rangle$ be the $p$-group of class $2$ described above. Then $\sigma$
defined by
$$a\mapsto ad^p, \hskip5mm b\mapsto b, \hskip5mm c\mapsto c,\hskip5mm d\mapsto
d$$
extends to a non-inner central automorphism of $H$. Let  $L$ be of order $p^3$
acting on $H$ via $\sigma$. Then $G$ is Miller if $N = \gen{[a,b]}$, and $G$ is
not if $N = \gen{[a,d]}$.

\begin{remark}
 The above construction of non-special Miller groups $G$ from special Miller
groups $H$ is valid even without hypotheses (i)-(iii) on $H$.  That
${\Aut}(G)={\Autcent}(G)$ can be proved using the  same arguments as in
\cite{yadav3}, which do not rely on hypotheses (i)-(iii). Then one can  apply
Theorem \ref{adney-yen} to show that $\Aut(G)$ is abelian.
\end{remark}

We conclude with the following remarks. In all known examples of Miller
$p$-groups $G$,  $G/Z(G)$ is homocyclic.
It will be interesting to know whether this happens in all Millers $p$-groups.
If true, one might expect if $\gamma_2(G)$ is always homocyclic. Again in the
known Millers $p$-groups $G$ one can observe that either $G'$ and $Z(G)$ have
same ranks or $G/Z(G)$ and $G/G'$ have same ranks. We wonder whether this is
true for
all Miller $p$-groups.

\bigskip
\newcommand{\Addresses}{{\footnotesize
Rahul Dattatraya Kitture,

\noindent  Post Doctoral Fellow, {\sc School of Mathematics, Harish-Chandra
Research
Institute, Chhatnag Road, Jhunsi, Allahabad 211019, INDIA \& Homi Bhabha
National Institute, Training School Complex, Anushakti Nagar, Mumbai 400085,
India} \textit{E-mail address: }  \texttt{rahul.kitture@gmail.com}

\vskip3mm
\noindent
Manoj K. Yadav,

\noindent  {\sc School of Mathematics, Harish-Chandra Research
Institute, Chhatnag Road, Jhunsi, Allahabad 211019, INDIA \& Homi Bhabha
National Institute, Training School Complex, Anushakti Nagar, Mumbai 400085,
India} \textit{E-mail address: } \texttt{myadav@hri.res.in}
}}

\noindent\Addresses
\end{document}